\documentclass[12pt]{article}

\usepackage{amssymb}
\usepackage{amsmath}
\usepackage{epsfig}

\newtheorem{thm}{Theorem}[section]
\newtheorem{lem}[thm]{Lemma}
\newtheorem{prp}[thm]{Proposition}
\newtheorem{cor}[thm]{Corollary}

\newenvironment{proof}{{\it Proof}.\ }{\hfill$\square$\par\medskip}

\def\cC{{\mathcal C}}
\def\cD{{\mathcal D}}

\def\frM{{\mathfrak M}}

\newcommand\NN{{\mathbb N}}

\newcommand\Sph{{\mathbb S}}
\newcommand\ZZ{{\mathbb Z}}

\def\to{\rightarrow}

\def\St{\operatorname{St}}

\newcommand\DMC[1]{\frM(#1)}
\newcommand\DMCp[1]{\frM_{\text{\rm pure}}(#1)}

\def\susp{\operatorname{susp}}
\def\hocolim{\operatorname{hocolim}}
\def\join{\star}

\title{Discrete Morse Complexes}
\author{Manoj K. Chari and Michael Joswig}
\date{\today.}

\begin{document}
\maketitle

\begin{abstract}
  We investigate properties of the set of discrete Morse functions on
  a simplicial complex as defined by Forman~\cite{896.57023}. It is
  not difficult to see that the pairings of discrete Morse functions
  of ~$\Delta$ again form a simplicial complex, the \emph{discrete
    Morse complex} of~$\Delta$.  It turns out that several known
  results from combinatorial topology and enumerative combinatorics,
  which previously seemed to be unrelated, can be re-interpreted in
  the setting of these discrete Morse complexes.  \medskip

\end{abstract}

\section{Introduction}
 
In the paper~\cite{896.57023}, Forman introduced the notion of a
discrete Morse function on an abstract simplicial complex and
developed a combinatorial analog of classical Morse theory. Discrete
Morse theory has proved to be an extremely useful tool in the study of
certain combinatorially defined spaces. In particular, it has been
applied quite successfully to the study of the topology of monotone
graph properties (see the papers of Babson et al.~\cite{936.57002},
Shareshian~\cite{Morse-2-connected} and
Jonsson~\cite{not-3-connected}). In this paper, we initiate the
investigation of the set of all possible Morse functions on a given
simplicial complex.  We give this set the structure of a simplicial
complex, which we call the {\it discrete Morse complex} associated
with the given simplicial complex.  This simplicial structure is
rather natural in the context of a graph-theoretical interpretation of
discrete Morse functions which is discussed in the aforementioned
papers of Chari~\cite{gshell}, Shareshian~\cite{Morse-2-connected} and
Jonsson~\cite{not-3-connected}.  For the special case of
one-dimensional simplicial complexes, that is graphs, discrete Morse
complexes reduce to complexes of rooted forests on graphs whose study
was undertaken independently by Kozlov~\cite{934.05041}.  In what
follows, we present our results in the topological and enumerative
study of discrete Morse complexes.

\section{Preliminaries}

Throughout the following let $\Delta$ be a finite abstract simplicial
complex. 

A {\em (discrete) Morse function\/} on~$\Delta$ is a function
$m:\Delta\to\NN$ with the following properties:  For each
$k$-face~$f\in\Delta$ there is at most one $(k+1)$-face~$g$
containing~$f$ with $m(g)<m(f)$, and there is at most one
$(k-1)$-face~$e$ contained in~$f$ with $m(e)>m(f)$.  The $k$-face~$f$
is {\em critical\/} with respect to~$m$ if $m$ attains a higher value
at all $(k+1)$-faces containing~$f$ and a lower value at all
$(k-1)$-faces contained in~$f$.  We can phrase it as follows: $f$ is
critical with respect to~$m$ if and only if, locally at~$f$, the
function~$m$ is strictly increasing with the dimension.

A key result of Forman is that if $m$ is a discrete Morse function
on~$\Delta$ with critical faces $f_1^{k_1},\ldots,f_n^{k^n}$, where
$\dim f_i^{k_i}=k_i$, then $\Delta$ is homotopic to a CW-complex with
$n$~cells of respective dimensions $k_1,\ldots,k_n$.  This is the
direct combinatorial equivalent to what is known from classical Morse
theory, for an introduction see Milnor's book~\cite{108.10401}.
Observe that the function $f\mapsto\dim f$ is a discrete Morse
function where all faces are critical.  In particular, a Morse
function always exists.  In order to understand the topological
structure of~$\Delta$ one needs a good Morse function, that is, a
Morse function with few critical faces.

For any  discrete Morse function $m$, it can be shown that for any non-critical face 
$f$, {\em exactly one} of the following is true:

 (i) there exists a (unique)  $(k+1)$-face~$g$
containing~$f$ with $m(g)<m(f)$,

 (ii) there exists a (unique) $(k-1)$-face~$g$
contained in~$f$ with $m(g) > m(f)$. 
  
Therefore, the set of non-critical faces with respect to~$m$ can be
uniquely partitioned in to pairs~$(f,g)$ where $f$ is a maximal face
of~$g$ and $m(f) > m(g)$.  Now, consider the Hasse diagram of~$\Delta$
as a directed graph; we direct all edges downward, that is, from the
larger faces to the smaller ones.  The previous observation implies
that the non-critical pairs form a matching in the Hasse diagram. If
we reverse the orientation of the arrows in this matching, it can be
shown that the resulting directed graph obtained is acyclic. We will
call such a matching on the Hasse diagram of the given simplicial
complex {\em acyclic}. Conversely, given an acyclic matching in the
Hasse diagram, one can construct a discrete Morse function with the
matching edges corresponding precisely to the non-critical pairs. The
critical faces are precisely those with no matching edges incident to
them.  We call two discrete Morse functions on~$\Delta$ {\em
  equivalent\/} if they induce the same acyclic matching.  In the
following we usually do not distinguish between equivalent discrete
Morse functions, that is, we identify a discrete Morse function with
its associated acyclic matching. For further details of this
interpretation of discrete Morse functions and its applications, we
refer to the papers of Chari~\cite{gshell},
Shareshian~\cite{Morse-2-connected} and
Jonsson~\cite{not-3-connected}.

The purpose of this paper is to study the set of all possible Morse
functions on a given simplicial complex by using this above
identification.  We define the {\em discrete Morse
  complex\/}~$\DMC\Delta$ of~$\Delta$, on the set of edges of the
Hasse diagram of~$\Delta$, as the set of subsets of edges of the Hasse
diagram of~$\Delta$ which form acyclic matchings. This is clearly an
abstract simplicial complex on the given vertex set.  Note that even
$\Delta$ is pure, the discrete Morse complex~$\DMC\Delta$ is not
necessarily pure itself.  Often it will be useful instead to consider
$\DMCp\Delta$, the {\em pure discrete Morse complex\/} of~$\Delta$,
the subcomplex of~$\DMC\Delta$ generated by the facets of maximal
dimension.  The facets of maximal dimension correspond to Morse
functions which are optimal in the sense that they lead to cell
decompositions with as few cells as possible.  Note that, for a
collapsible simplicial complex such an optimal Morse function
corresponds to a collapsing strategy (up to a reordering of the
elementary collapses) and vice versa.

\section{Discrete Morse Complexes of Graphs}

As can be expected, discrete Morse complexes are typically very large
and very complicated spaces. To obtain some sort of intuition about
these spaces, it is helpful to consider the one-dimensional case, that
is, graphs. We first observe that the (undirected) Hasse diagram of a
graph~$\Gamma$ is obtained by subdividing each edge of the graph exactly once.
Now a Morse matching on such a complex gives us pairs (of non-critical
faces) which are all of the type $(v,e)$ where $v$ is a node in $\Gamma$
and $e$ is an edge of $\Gamma$ with $v$ as one of its end points. Consider
the subgraph $S(M)$ of $\Gamma$ all such edges $e$ (with both endpoints
included) which appear in $M$ and orient each edge $e$ away from $v$
in $S(M)$. The matching property applied to $e$ ensures that this
construction is well defined while the matching property at each node
$v$ ensures that the out-degree at each node is at most one. From the
acyclic property of $M$, we can deduce that the subgraph $F(M)$ as an
undirected graph contains no cycles and hence is a forest. Since the
outdegree at each node of $S(M)$ is at most one, each component of the
forest has a unique ``sink''(often called a ``root'') with respect to
the given orientation. Given any graph $\Gamma$, we call an oriented subset
$F$ of edges of $G$, a {\it rooted forest} of $\Gamma$ if $F$ is forest as
an undirected graph and further, every component of $F$ has a unique
root with respect to the given orientation.  We have argued above that
every acyclic  matching for the Hasse diagram of a graph corresponds in a natural way to a
rooted forest in a graph and this can be easily reversed to yield the following.

\begin{prp}
 The set of Morse functions on a graph $\Gamma$ is in one-to-one
 correspondence with the set of rooted forests of $\Gamma$.
\end{prp}

Complexes of rooted trees and forests of graphs have been
independently investigated by Kozlov~\cite{934.05041} and the above
proposition shows that the complexes he considers are, in fact,
discrete Morse complexes of graphs.
 
In particular, the facets of the complex $\DMC{\Gamma}$ correspond
precisely to the rooted spanning trees of $\Gamma$. Each facet gives
rise to a ``different'' proof of the elementary fact that a
(connected) graph with $m$ edges and $n$ nodes is homotopy equivalent
to a wedge of $m-n+1$ circles. The rooted spanning tree (in isolation)
can be collapsed according to the orientation to a point represented
by the root node. The rest of the $m-n+1$ edges form the $m-n+1$
critical $1$-cells giving the homotopy type.  This simple result
yields some interesting enumerative consequences for the $f$-vector of
the complex $\DMC{\Gamma}$ which we now discuss.  Recall that the
$f$-vector of a complex just lists for all $i$, the number~$f_{i}$ of
$i$-dimensional faces.  From well-known formulae for the number of
rooted forests on $n$ nodes with $k$-trees, we get an explicit formula
for the $f$-vector of $\DMC{K_{n}}$.

\begin{cor} 
  The $f$-vector of the discrete Morse complex of the complete graph
  $K_{n}$ on~$n$ nodes is given by
  $$ f_{i-1} = \binom{n}{i} (n-i) n^{i-1}. $$
\end{cor} 

For general graphs, we can relate the $f$-vector of $\DMC{\Gamma}$ to
the characteristic polynomial of the Laplacian matrix of the graph
using the proposition above. The spectrum of the Laplacian is a
fundamental algebraic object associated with a graph and has been
studied extensively (see Biggs~\cite{797.05032}). In what follows, we
assume familiarity with the basic notions of algebraic graph theory
and we will follow the notation and terminology of
Biggs~\cite{797.05032}. Given a (connected) graph $\Gamma$ (which we
assume for convenience to be connected), let ${\bf Q}(\Gamma)$ be the
Laplacian of the graph and let $\sigma(\Gamma;\mu)$ be the
characteristic polynomial of ${\bf Q}(\Gamma)$ given by
$\sigma(\Gamma;\mu) = det(\mu{\bf I}-{\bf Q}(\Gamma))$.

\begin{cor}
  If $(f_0,f_1,\ldots)$ is the $f$-vector of $\DMC{\Gamma}$ then
  $$\sigma(\Gamma;\mu) = \sum f_{i-1} (-1)^i \mu^{n-i}.$$
\end{cor}
The proof of this is immediate from the above proposition and Theorem
7.5 of \cite{797.05032}

Now we move on the topological properties of discrete Morse complexes
of graphs.  We will frequently use a certain special case of a general
result about homotopy colimits of diagrams of spaces, cf.\ Welker,
Ziegler, \v{Z}ivaljevi\'c~\cite{hocolim}.

\begin{prp}\label{lem:union}
  Let $A$, $B$ be subcomplexes of~$\Delta$ such that both inclusion
  maps $A\cap B\hookrightarrow A$ and $A\cap B\hookrightarrow B$ are
  homotopic to the constant map.

  Then $A\cup B\simeq A\vee\susp(A\cap B)\vee B$.
\end{prp}

\begin{proof}
  Consider the diagram~$\cD$ associated to the poset $(\{A\cap
  B,A,B\},\ge)$.  Note that inclusion is reversed such that $A\cap B$
  becomes~$\widehat{1}$.  The inclusion maps $A\cap B\hookrightarrow
  A$ and $A\cap B\hookrightarrow B$ clearly are closed cofibrations.
  By the Projection Lemma~\cite[4.5]{hocolim} we have that $A\cup
  B\simeq\hocolim\cD$.
  
  The inclusion maps being homotopic to the constant map, it follows
  from the Wedge Lemma \cite[4.9]{hocolim} that $\hocolim\cD\ \simeq\ 
  A*\emptyset\,\vee\,(A\cap B)*\Sph^0\,\vee\,B*\emptyset\ \simeq\ 
  A\vee\susp(A\cap B)\vee B$.
\end{proof}

Two instances of the preceding proposition are particularly relevant
for our discussion.

\begin{cor}\label{cor:union-contractible}
  Let $A$, $B$ be contractible.

  Then $A\cup B\simeq\susp(A\cap B)$.
\end{cor}

\begin{cor}\label{cor:union-spheres}
  Assume that $A\simeq\Sph^n\simeq B$ and $A\cap B\simeq\Sph^r$ with $r<n$.
  
  Then $A\cup B\simeq\Sph^{r+1}\vee\Sph^n\vee\Sph^n$.
\end{cor}

\section{The Pure Discrete Morse Complex of a Circle}

\begin{figure}[htbp]
  \begin{center}
    \input{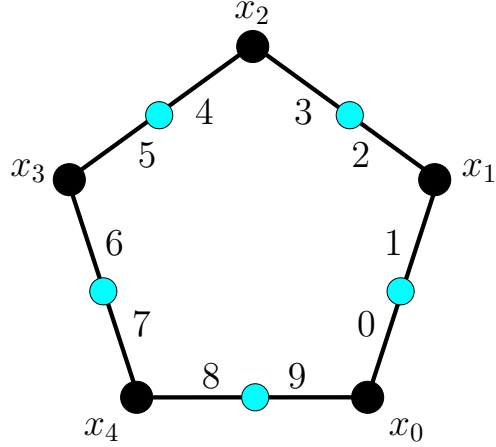}
    \caption{Graph~$C_5$ and numbering of the vertices of~$\DMC{C_5}$.}
    \label{fig:circle}
  \end{center}
\end{figure}

Let $C_n$ be the cyclic graph on $n$~nodes.  Obviously, $C_n$ is
homeomorphic to the circle~$\Sph^1$.  We choose the following
notation.  The nodes of~$C_n$ are denoted by $x_0,x_1,\ldots,x_{n-1}$
with edges $(x_i,x_{i+1})$; all indices are taken modulo~$n$.  The
$2n$~vertices of~$\DMC{C_n}$ are identified with numbers
$0,1,2,\ldots,2n-1$ such that the vertex~$2i$ corresponds to the
pair~$(x_i,(x_i,x_{i+1}))$ and $2i+1$ corresponds to the
pair~$(x_{i+1},(x_i,x_{i+1}))$, see Figure~\ref{fig:circle}.

The discrete Morse complex~$\DMC{C_n}$ is pure if and only if $n\le5$.
Its $f$-vector is known to be
$$f_i=\frac{2n}{2n-i-1}\binom{2n-i-1}{i-1},$$
see
Stanley~\cite[2.3.4]{EnumComb1}. We note that for any $n$, the pure
discrete Morse complex~$\DMCp{C_n}$ is of dimension $n-2$.

\begin{thm}\label{thm:circle}
  Let $n\ge 4$.
  
  Then the pure discrete Morse complex~$\DMCp{C_n}$ is homotopic to
  $\Sph^2\vee\Sph^{n-2}\vee\Sph^{n-2}$.
\end{thm}

Kozlov computed the homotopy type of the discrete Morse
complex~$\DMC{C_n}$.  Our result~\ref{thm:circle} was obtained
independently. 

\begin{thm}{\rm (Kozlov~\cite[Proposition 5.2]{934.05041})}
\label{thm:Kozlov}
  $$\DMC{C_n}\simeq\left\{\begin{array}{cl}
      \Sph^{2k-1}\vee\Sph^{2k-1}\vee\Sph^{3k-2}\vee\Sph^{3k-2},
      &\text{if $n=3k$;}\\
      \Sph^{2k}\vee\Sph^{3k-1}\vee\Sph^{3k-1},
      &\text{if $n=3k+1$;}\\
      \Sph^{2k}\vee\Sph^{3k}\vee\Sph^{3k},
      &\text{if $n=3k+2$.}\end{array}\right.$$
\end{thm}

Before we prove Theorem~\ref{thm:circle}, we will establish the
following useful lemma.

\begin{lem}
  The pure Morse complex of any path (with more than two nodes) is
  collapsible.
\end{lem}

\begin{proof}
  The poset of faces of a path is obtained by a subdividing of each
  edge of the path, when constructing Morse matchings, the acyclic
  property is trivially satisfied.  Therefore, the pure Morse complex
  of the path with $n$ edges is simply the complex of partial
  matchings that are extendible to perfect matchings for the path with
  $2n$ edges.  We will call this the {\it pure matching complex} for
  the path with $2n$ edges.  Assuming that the $2n$ edges are labeled
  $1,2,\ldots, 2n$, it is clear that there is exactly one perfect
  matching which contains $2$ , namely $ F = \{2,4,\ldots, 2n\}$ and
  every other perfect matching contains the edge $1$. It is easy to
  see that $F$ can be collapsed onto the face $\{4,\ldots, 2n\}$ which
  is contained in the facet $\{1,4,\ldots, 2n\}$. As a result, we are
  left with a cone with apex~$1$, which is obviously collapsible.
  
  Note that a path with exactly two nodes, that is, an interval has a
  discrete Morse complex isomorphic to~$\Sph^0$.
\end{proof}

\begin{proof} (of Theorem~\ref{thm:circle}) 
  Observe that $\DMCp{C_n}=\St0\cup\St1\cup\St2\cup\St3$.
  
  We claim that $\St0\cap\St2=\St02\cup\{4,6,8,\ldots,2n-2\}$. To show
  this assume that $F$ is a maximal face of $\St0\cap\St2\setminus\St02$.
  Then both $F\cup 0$ and $F\cup 2$ are facets of $\St0$ and $\St2$,
  respectively. It follows that $F$ consists of $n-2$ elements and it
  cannot possibly contain $1$, $3$ or~$2n-1$ and therefore these
  elements must come from $\{4,5,6,\ldots,2n-2\}$.  From the matching
  property that is required it follows that there exactly one
  possibility is, that is, $F = \{4,6,8,\ldots,2n-2\}$. Also, to any
  proper subset of $F$ , say $G$, one can always add the edges~$0$
  and~$2$ to obtain a Morse matching, which is obviously an element
  of~$\St02$. This completes the proof of the claim.  Now the
  star~$\St02$ is clearly contractible, and as shown above the
  boundary of the face~$\{4,6,8,\ldots,2n-2\}$ is entirely contained
  in~$\St02$.  We infer that $\St0\cap\St2$ is homotopic
  to~$\Sph^{n-3}$.  Similarly, $\St1\cap\St3\simeq\Sph^{n-3}$.
  
  As the stars~$\St0$ and~$\St2$ are contractible, we can apply
  Corollary~\ref{cor:union-contractible} to derive that $\St0\cup\St2$
  is homotopic to the suspension of the intersection
  $\St0\cap\St2\simeq\Sph^{n-3}$.  Thus
  $\St0\cup\St2\simeq\Sph^{n-2}$.  Similarly,
  $\St1\cup\St3\simeq\Sph^{n-2}$.
 
  Now we consider $(\St0\cup\St2)\cap(\St1\cup\St3)$ = $(\St0\cap\St1)
  \cup (\St0\cap\St3) \cup (\St2\cap\St1) \cup (\St2\cap\St3)$, which
  we claim is equal to $\St03\cup B$, where $B$ is the pure matching
  complex of the path with edges $\{4,5,\ldots,2n-1\}$. To show this,
  assume $F$ is a facet of the intersection. It is easy to see from
  the matching requirement that only way that $F$ can have $n-1$
  elements is if $ 0 \in F$ and $3 \in F$, that is $F$ is a facet of
  $\St03$. It is also evident that in this instance, the set $G = F\setminus
  \{0,3 \}$ is a facet of the pure matching complex of the path on the
  edges $\{5,6,\ldots,2n-2\}$ .
  
  Now the facets of $(\St0\cup\St2)\cap(\St1\cup\St3)\setminus\St03$ are
  subsets of $\{4,5,\ldots,2n-1\}$).  It is clear that any such facet
  $F$ is also a facet of the pure matching complex of the path on the
  edges $\{4,5,\ldots,2n-1\}$. Conversely, to any facet $F$ of the pure
  matching complex of this path we can add either $0$ or $2$ to $F$ to
  get a facet of $(\St0\cup\St2)$ {\it and} we can add either $1$ or
  $3$ to get a facet of $(\St1\cup\St3)$ so that $F$ will be facet of
  $(\St0\cup\St2)\cap(\St1\cup\St3)$.
 
  Thus we have $(\St0\cup\St2)\cap(\St1\cup\St3)$ is the union of the
  two contractible complexes $\St03$ and $B$.  Note that if we have a
  facet of the pure matching complex of the path on $2n$~edges,
  consecutively labeled starting with an odd number, then every facet
  consists of a string of odd edges (possibly empty) followed by an
  even string of vertices (possibly empty). On the other hand, if the
  labeling starts with an even number, then the facets consist of an
  even string followed by an odd string. This observation is useful in
  determining $\St03\cap B$, which following the above arguments, is
  the intersection of the pure matching complex on the path
  $\{4,5,\ldots,2n-1\}$ with that of the pure matching complex on the
  path $\{5,\ldots,2n-2\}$.  It follows that any face in this
  intersection consists entirely of odd vertices or entirely of even
  vertices.  Now the two faces $\{6,8,\ldots,2n-2\}$ and
  $\{5,7,\ldots,2n-1\}$ are in the intersection and obviously, they
  are the unique maximal even and odd sets, respectively, in the
  intersection.  Therefore, we have shown intersection $\St03\cap B$
  is a disjoint union of two non-empty simplices.  That is, it is
  homotopic to~$\Sph^0$.  Due to
  Corollary~\ref{cor:union-contractible} we have that $\St03\cup
  B=(\St0\cup\St2)\cap(\St1\cup\St3)\simeq\Sph^1$.

  Finally, the claim follows from Corollary~\ref{cor:union-spheres}.
\end{proof}

We remark here the small cases of $n$  can also be treated by shelling techniques.

\section{The Discrete Morse Complex of the Simplex}

Let $\Delta_d$ be the $d$-dimensional simplex, that is $\Delta_d$ is
the Boolean lattice on $d+1$~points.  The Hasse diagram of~$\Delta_d$
is isomorphic to the graph of the $(d+1)$-dimensional cube: As the
vertices of the cube take all $0/1$-vectors of length~$d+1$; the
linear function $x_0+\dots+x_d$ induces an acyclic orientation on the
graph of the $0/1$-cube.  Mapping a subset of $\{0,1,\ldots,d\}$ to
its characteristic function yields the desired isomorphism.  By
$\Gamma_{d+1}$ we will denote the directed graph of the $(n+1)$-cube
whose arcs point towards lower values of the named linear function.

Except for very small~$d$ it seems to be extremely difficult to
determine the topological types of~$\DMC{\Delta_d}$
and~$\DMCp{\Delta_d}$.  It even seems to be hard to compute the
generating function of the $f$-vector.

\begin{prp}
  The discrete Morse complex of~$\Delta_d$ is homotopic to
  \begin{enumerate}
  \item[a)] the $0$-sphere~$\Sph^0$ if $d=1$,
  \item[b)] the wedge $\Sph^1\vee\Sph^1\vee\Sph^1\vee\Sph^1$ if $d=2$.
  \end{enumerate}
\end{prp}

\begin{proof}
  A line segment can be oriented in two different ways.  Hence the
  result for $d=1$.

  \begin{figure}[htbp]
    \begin{center}
      \epsfig{file=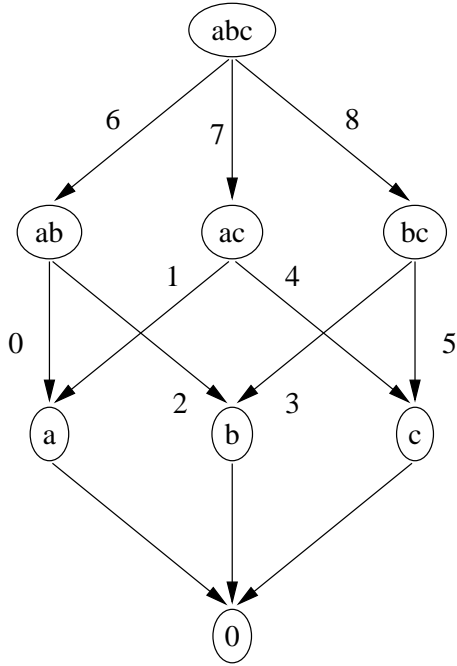,width=6cm}
      \caption{Hasse diagram of the $2$-simplex.}
    \end{center}
    \label{fig:2-simplex}
  \end{figure}

  Now consider the $2$-simplex with vertices $a$, $b$, $c$; label the
  edges of the Hasse diagram from $0$ to~$8$ as in
  Figure~\ref{fig:2-simplex}.  For any subset $X\subseteq\Delta_2$ let 
  $\frM[X]$ be the discrete Morse complex of the subcomplex generated
  by~$X$.  Moreover, $\frM=\frM[\Delta_2]$.

  Clearly, each maximal Morse matching contains either $6$ or~$7$
  or~$8$, that is, $\frM=\St 6\cup\St 7\cup\St 8$.  Now $\St
  6=6\join\frM[ac,bc]$, $\St 7=7\join\frM[ab,bc]$, and $\St
  8=8\join\frM[ab,ac]$.
  
  Moreover, the intersection $\St 6\cap\St 7$ equals $\frM[bc]$ which
  consists of two isolated points.  Thus $\St 6\cup\St 7\simeq\Sph^1$.
  
  Observe that $(\St 6\cup\St 7)\cap\St 8=\frM[ab]\cup\frM[ac]$ has
  four isolated points, that is, it is equal
  to~$\Sph^0\stackrel\cdot\cup\Sph^0$.  We infer that
  $\frM\simeq\Sph^1\vee\susp(\Sph^0\stackrel\cdot\cup\Sph^0)
  \simeq\bigvee_4\Sph^1$.
\end{proof}

\begin{figure}[htbp]
  \begin{center}
    \epsfig{file=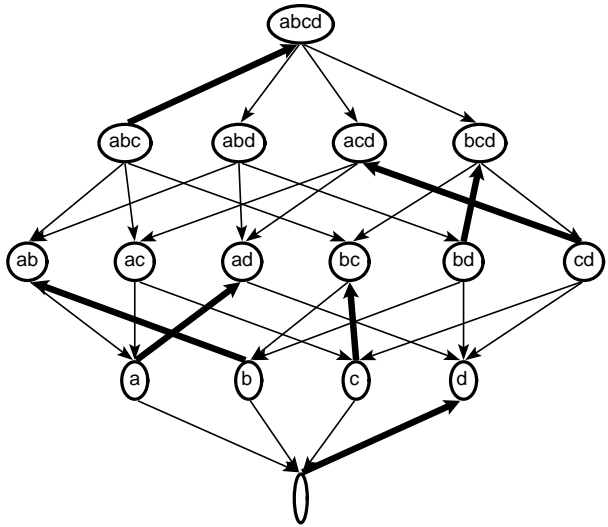,width=8cm}
    \caption{Maximal acyclic matching in the Hasse diagram of the $3$-simplex.}
    \label{fig:3-simplex}
  \end{center}
\end{figure}

Note that for $n\ge3$ the Morse complex of~$\Delta_n$ is no longer
pure.  For an example of a maximal Morse matching of the $3$-simplex
which is not perfect see Figure~\ref{fig:3-simplex}.

We computed the homology of the (pure) discrete Morse complex of the
$3$-simplex using two independent software implementations by
Heckenbach~\cite{homology} and by Gawrilow and
Joswig~\cite{SimplicialHomology}.  See Table~\ref{tab:3-simplex} for
the results.

\begin{table}
  \begin{center}
    \begin{tabular}{l|cc}
      & $f$-vector & reduced integer homology\\
      \hline
      $\DMC{\Delta_3}$  & $(28,300,1544,3932,4632,2128,256)$ & $(0,0,0,0,\ZZ^{99},0,0)$\\ 
      $\DMCp{\Delta_3}$ & $(28,300,1544,3680,3672,1600,256)$ & $(0,0,0,\ZZ^{81},0,0,0)$
    \end{tabular}
    \caption{Discrete Morse complex of the $3$-simplex.  All vectors
      are written from left to right with increasing dimension.}
    \label{tab:3-simplex}
  \end{center}
\end{table}

We now turn our attention to the number of faces of discrete Morse
complex of the simplex. For collapsible complexes, such as the
simplex, the \emph{perfect Morse matchings}, that is, the facets of
the pure part of the discrete Morse complex, correspond to collapsing
strategies (modulo ordering of the elementary collapses).  Given a
perfect Morse matching~$\mu$ of any collapsible complex each $k$-face
is paired to a unique $l$-face where $l\in\{k-1,k+1\}$.  Now let $\mu$
be a perfect Morse matching of the $(n-1)$-simplex.  Define $T_k(\mu)$
to be the set of $k$-faces of~$\Delta_{n-1}$ which are paired to
$(k-1)$-faces (via~$\mu$). We now relate these sets to certain
subcomplexes to simplices which were studied by Kalai~\cite{535.57011}.

A \emph{$(k,n)$-tree~$T$} is a subset of~$\Delta_{n-1}^k$ of
cardinality~$\binom{n-1}{k}$ which has the property that for
$\Delta(T)=\Delta_{n-1}^{\le k-1}\cup T$ we have $H_k(\Delta(T))=0$.
From an Euler characteristic argument we infer that that
$H_{k-1}(\Delta(T))$ is finite. The following lemma is immediate from
the definition of $T_k(\mu)$.
\begin{lem}
  For each $\mu$ and each~$k$ the set $T_k(\mu)$ is a $(k,n)$-tree.
\end{lem}

\begin{lem}\label{lem:1-1}
  Let $\mu$, $\mu'$ be perfect Morse matchings of~$\Delta_{n-1}$ with
  $(T_{n-1}(\mu),\dots,T_0(\mu))=(T_{n-1}(\mu'),\dots,T_0(\mu'))$.

  Then $\mu=\mu'$.
\end{lem}

\begin{proof}
  On the contrary assume that $\mu\ne\mu'$.  Hence, for some~$k$ the
  restrictions $\left.\mu\right|_{T_k(\mu)}$ and
  $\left.\mu'\right|_{T_k(\mu')}$ differ.  Abbreviate $T_k$ for
  $T_k(\mu)=T_k(\mu')$ and $T_{k-1}$ for $T_{k-1}(\mu)=T_{k-1}(\mu')$,
  respectively.  Let $\Gamma$ be the subgraph of the Hasse diagram
  of~$\Delta_{n-1}$ which is induced on the vertex set~$T_k\cup
  (\Delta^{k-1}\setminus T_{k-1})$.  Both, $\mu$ and~$\mu'$ induce
  perfect matchings of the bipartite graph~$\Gamma$.  So their
  symmetric difference is a union of cycles.  We arrive at a
  contradiction because of the acyclicity condition on Morse
  matchings.
\end{proof}

Let $\cC(n,k)$ be the set of $(k,n)$-trees.

\begin{thm}{\rm (Kalai~\cite{535.57011})}\label{thm:Kalai}
  For arbitrary $n$ and~$k$ we have
  $$\sum_{C\in\cC(n,k)}\left|H_{k-1}(\Delta(C))\right|^2=n^{\binom{n-2}{k}}$$
  and

  $$\left|\cC(n,k)\right|\le\left(\frac{en}{k+1}\right)^{\binom{n-1}{k}},$$
  where $e$ is Euler's constant.
\end{thm}

The preceding results immediately yield an upper bound on the
number~$f(n)$ of perfect Morse matchings of the $n$-simplex.

\begin{cor}\label{cor:MorseMatchingsSimplex}
  The number of perfect Morse matchings of the $n$-simplex is bounded
  from above by $$f(n)\le(n+1)^{2^{n-1}}.$$
\end{cor}

Note that, $f(1)=2$, $f(2)=9$, and $f(3)=256$; that is, in principal,
the upper bound is tight.  However, for larger values of~$n$ the
estimate becomes increasingly inaccurate for two obvious reasons.
Firstly, the formula in Theorem~\ref{thm:Kalai} also counts
$(k,n)$-trees~$T$ for which~$\Delta(T)$ is not collapsible.  Secondly,
each summand is weighted whereas here we are only interested in the
number of summands.

So it seems reasonable to look for a better upper bound. A possible
way is straightforward from the definition of a perfect Morse matching
as a special type of perfect matching.  This leads to the problem of
counting perfect matchings in the graph of the $(n+1)$-dimensional
cube.  There is an asymptotic solution to this problem, which is due
to Clark, George, and Porter~\cite{901.05056}. Here we are interested
only in the upper bound.

\begin{thm}{\rm (Clark, George, and Porter)}
  \label{thm:matchings-cube} The number of perfect matchings
  of the graph of the $(n+1)$-dimensional cube is bounded from above
  by $$f(n)\le{(n+1)!}^{\frac{2^n}{n+1}}.$$
\end{thm}

Unfortunately, a direct computation shows that the bound from
Corollary~\ref{cor:MorseMatchingsSimplex} is always better than the
one derived from Theorem~\ref{thm:matchings-cube}

What about lower bounds?  The interpretation of the Hasse diagram of
the $n$-simplex as the directed graph~$\Gamma_{n+1}$ of the
$(n+1)$-cube suggests a way of constructing perfect Morse matchings
recursively.  Recall that the vertices of~$\Gamma_{n+1}$ are the
vertices of the $(n+1)$-dimensional $0/1$-cube.  For arbitrary
$i\in\{0,\ldots,n\}$ the vertices satisfying the equation $x_i=0$ are
the vertices of an $n$-cube, whose graph we denote
by~$\Gamma_{n+1,i}^-$; similarly, we obtain the
graph~$\Gamma_{n+1,i}^+$ of another $n$-cube for $x_i=1$.  We call
$\Gamma_{n+1,i}^+$ and $\Gamma_{n+1,i}^-$ \emph{bottom} and
\emph{top}, respectively.  Observe that all arcs in between point from
bottom to top.  Thus any perfect acyclic matching
of~$\Gamma_{n+1,i}^+$, combined with any perfect acyclic matching
of~$\Gamma_{n+1,i}^-$ yields a perfect acyclic matching
of~$\Gamma_{n+1}$.  In principal, we can do this for
every~$i\in\{0,\ldots,n\}$, but we may obtain the same matching for
different~$i$.

\begin{prp}\label{prp:lower-bound}
  Let $r(1)=1$, $r(2)=2$, $r(3)=9$, and, for $n\ge3$, recursively, 
  $$r(n+1)=\frac{(n+1)(n-1)}{n}\ r(n)^2.$$
  
  The number of perfect Morse matchings of the $n$-simplex is bounded
  from below by $$f(n)\ge r(n+1).$$
\end{prp}

\begin{proof}
  The numbers of perfect matchings of the graph of the $(n+1)$-cubes
  for $n+1\le3$ are easy to determine.  All these matchings are
  acyclic and thus are Morse matchings of the respective $n$-simplex.
  
  We say that an edge of~$\Gamma_{n+1}$ is \emph{in direction}~$i$ if
  its vertices differ in the $i$-th coordinate.  In the following we
  construct perfect acyclic matchings of cubes which contain edges of
  all but one direction.  Observe that all perfect matchings of
  $\Gamma_1$, $\Gamma_2$ and~$\Gamma_3$ are of this kind.
  
  Choose $k\in\{0,\ldots,n\}$.  Suppose we have two such
  matchings~$\mu^+$, $\mu^-$ in~$\Gamma_{n+1,k}^+$ and
  $\Gamma_{n+1,k}^-$, respectively.  Then $\mu=\mu^+\cup\mu^-$ is a
  perfect acyclic matching of~$\Gamma_{n+1}$.  Now $\mu$ contains edges
  of either $n-1$ or~$n$ directions.  Fix $\mu^+$, and let
  $i\in\{0,\ldots,k-1,k+1,\ldots,n\}$ be the unique direction which
  $\mu^+$ does not contain an edge of.  Now $\mu$ contains edges from
  $n$~directions if and only if $\mu^-$ contains edges of
  direction~$i$.
  
  If there are $r$ perfect acyclic matchings of the $n$-cube
  containing edges of all but one direction, then $r(n-1)/n$ of them
  contain edges of a given direction.  This gives $r^2(n-1)/n$
  different perfect acyclic matchings of the $\Gamma_{n+1}$, which
  contain edges from all but direction~$n$.
  
  Now there are $n+1$ choices for~$k$, and all of them yield different
  matchings.
\end{proof}

The number $p(n)$ of perfect matchings of the graph of the $n$-cube is
known for small values of~$n$: In addition to the obvious values
Graham and Harary~\cite{626.05043} computed $p(4)=272$ and
$p(5)=589,185$.  Moreover, in~\cite{901.05056} the authors mention
that Weidemann showed that $p(6)=16,332,454,526,976$.

Note that the graph of the $4$-cube has a perfect acyclic matching
using all directions; it can be constructed from a Hamiltonian cycle.

In order to give a vague idea about the growth of the function $r$
from Proposition~\ref{prp:lower-bound} one can unroll the recursion.

\begin{cor}
  The number of perfect Morse matchings of the $n$-simplex is bounded
  from below by $$f(n)\ge r(n+1)>\prod_{k=1}^{n-1} k^{2^{n-k-1}}.$$
\end{cor}

\begin{proof}
  We will prove the result by induction on~$n$. The initial case
  $1^{2^0}=1<2=r(2)$ is clear.  Further,

  \begin{eqnarray*}
    \prod_{k=1}^{n-1} k^{2^{n-k-1}}&=&(n-1)\left[\prod_{k=1}^{n-2}
      k^{2^{n-k-2}}\right]^2\\
    &<&(n-1)\ r(n)^2\\
    &<&\frac{(n+1)(n-1)}{n}\ r(n)^2=r(n+1).
  \end{eqnarray*}
\end{proof} 

This way we obtain a growth rate for the number of perfect Morse
matchings of the $n$-simplex which is approximately $(1.289)^{2^n}$.
We conjecture that the precise value of $f(n)$ has a function of~$n$
which goes to infinity with~$n$ as the base of this double exponent.

\section*{Acknowledgements}

We owe much to Volkmar Welker for stimulating discussions leading to
the present paper.  We are also indebted to Volker Kaibel for valuable
suggestions concerning the construction of matchings in cube graphs.

\bibliographystyle{amsplain}
\bibliography{math,top,graph,poly,morse,soft,mic}
\vspace{1cm}

\begin{minipage}[t]{9cm}
  Manoj K. Chari\\
  Dept.\ of Mathematics\\
  Louisiana State University\\
  Baton Rouge, LA 70803-4918\\
  U.S.A.
\end{minipage}
\begin{minipage}[t]{9cm}
  Michael Joswig\\
  Fachbereich Mathematik, MA 7-1\\
  Technische Universit\"at Berlin\\
  Stra\ss{}e des 17. Juni 136\\
  10623 Berlin\\
  Germany
\end{minipage}
\end{document}